\pdfoutput=1 

\documentclass[11pt]{amsart}

\usepackage{amssymb,amsthm,enumitem,colonequals,tikz-cd,microtype,quiver,mathtools,graphicx}
\usepackage[normalem]{ulem} 

\usepackage[osf]{Baskervaldx}
\usepackage[baskervaldx]{newtxmath}
\usepackage[cal=boondoxo]{mathalfa} 

\usepackage[top=3.75cm, bottom=3cm, left=3.5cm, right=3.5cm]{geometry}

\usepackage{xcolor}
\colorlet{darkblue}{blue!55!black}
\colorlet{darkcyan}{cyan!50!black}
\colorlet{darkgreen}{green!60!black}

\PassOptionsToPackage{hyphens}{url}
\usepackage{hyperref}
\hypersetup{
    colorlinks=true,
    linkcolor=darkblue,
    urlcolor=darkcyan,
    citecolor=darkgreen,
}

\def\eqref#1{\textcolor{darkblue}{(\ref{#1})}}

\usepackage[nameinlink]{cleveref} 
\Crefformat{section}{#2\S#1#3}
\Crefmultiformat{section}{#2\S\S#1#3}{ and~#2#1#3}{, #2#1#3}{, and~#2#1#3}

\usepackage[pagewise]{lineno}
\let\oldequation\equation
\let\oldendequation\endequation
\renewenvironment{equation}{\linenomathNonumbers\oldequation}{\oldendequation\endlinenomath}
\expandafter\let\expandafter\oldequationstar\csname equation*\endcsname
\expandafter\let\expandafter\oldendequationstar\csname endequation*\endcsname
\renewenvironment{equation*}{\linenomathNonumbers\oldequationstar}{\oldendequationstar\endlinenomath}
\let\oldalign\align
\let\oldendalign\endalign
\renewenvironment{align}{\linenomathNonumbers\oldalign}{\oldendalign\endlinenomath}
\expandafter\let\expandafter\oldalignstar\csname align*\endcsname
\expandafter\let\expandafter\oldendalignstar\csname endalign*\endcsname
\renewenvironment{align*}{\linenomathNonumbers\oldalignstar}{\oldendalignstar\endlinenomath}

\theoremstyle{plain}
\newtheorem{theorem}{Theorem}[section]
\newtheorem{lemma}[theorem]{Lemma}

\newtheorem{proposition}[theorem]{Proposition}

\theoremstyle{definition}
\newtheorem{definition}[theorem]{Definition}
\newtheorem{example}[theorem]{Example}
\newtheorem{remark}[theorem]{Remark}

\newtheorem*{ack}{Acknowledgments}

\AddToHook{env/lemma/begin}{\crefalias{theorem}{lemma}}
\AddToHook{env/corollary/begin}{\crefalias{theorem}{corollary}}
\AddToHook{env/proposition/begin}{\crefalias{theorem}{proposition}}
\AddToHook{env/definition/begin}{\crefalias{theorem}{definition}}
\AddToHook{env/remark/begin}{\crefalias{theorem}{remark}}
\AddToHook{env/example/begin}{\crefalias{theorem}{example}}

\numberwithin{equation}{section}
\numberwithin{theorem}{section}

\title[Remarks on Fourier--Mukai invariances]{Remarks on Fourier--Mukai invariances}

\author[E.~Guisado Villalgordo]{El\'{i}as Guisado Villalgordo}
\address{E.~Guisado Villalgordo,
BCAM -- Basque Center for Applied Mathematics, Mazarredo 14, 48009
Bilbao, Basque Country -- Spain}
\email{eguisado@bcamath.org}

\author[P.~Lank]{Pat Lank}
\address{P.~Lank,
Dipartimento di Matematica “F. Enriques”, Universit\`{a} degli Studi di Milano, Via Cesare
Saldini 50, 20133 Milano, Italy}
\email{lankp.math@gmail.com}

\date{\today}



\begin{document}
    
\begin{abstract}
    This note proves that smoothness and Gorensteinness are Fourier--Mukai invariants for proper flat families of Noetherian algebraic spaces.
    Moreover, a kernel realizing such a partnership is perfect if, and only if, the families are smooth.
\end{abstract}

\maketitle

\tableofcontents

\section{Introduction}
\label{sec:intro}

Fourier--Mukai partnerships provide rich comparisons of geometric objects through their derived categories. 
It first appeared in the work of \cite{Mukai:1981}, and has been extended to relative settings \cite{Ruiperez/Hernandez/Martin/SanchodeSalas:2009,GuisadoVillaalgordo/Lank:2026,GuisadoVillalgordo/Lank/ManaliRahul/Pavic:2025}. 
A fundamental question is which geometric properties are preserved under Fourier--Mukai partnerships.  

We first fix notation.
Let $f_i\colon Y_i\to S$ be proper flat morphisms of Noetherian algebraic spaces for $i=1,2$.
Denote by $p_i\colon Y_1\times_S Y_2\to Y_i$ the projections.
For $E\in D^b_{\operatorname{coh}}(Y_1\times_S Y_2)$, the associated \textit{integral transform} is $\Phi_E\colon D_{\operatorname{qc}}(Y_1)\to D_{\operatorname{qc}}(Y_2)$ where $A\mapsto \mathbf{R}(p_2)_\ast(\mathbf{L}(p_1)^\ast A\otimes^{\mathbf{L}}E)$. 
We say that $Y_1$ and $Y_2$ are \textit{Fourier--Mukai $S$-partners} if $\Phi_E$ restricts to an equivalence $D^b_{\operatorname{coh}}(Y_1)\to D^b_{\operatorname{coh}}(Y_2)$ for some such kernel $E\in D^b_{\operatorname{coh}}(Y_1\times_S Y_2)$.

\subsection{Results}
\label{sec:intro_results}

For smooth projective schemes over a field, kernels of such transforms are perfect complexes.
For nonsmooth families, perfectness need not hold.

\begin{theorem}
	\label{thm:smoothness}
	Let $f_i\colon Y_i\to S$ be proper flat morphisms of Noetherian algebraic spaces.
	Suppose $Y_1$ and $Y_2$ are Fourier--Mukai $S$-partners given by a kernel $E\in D^b_{\operatorname{coh}}(Y_1\times_S Y_2)$.
	Then the following are equivalent:
	\begin{enumerate}
		\item $f_1$ is smooth 
		\item $f_2$ is smooth 
		\item $E\in \operatorname{Perf}(Y_1\times_S Y_2)$.
	\end{enumerate}
\end{theorem}

\begin{proof}
	This follows from \Cref{prop:smoothness,prop:kernel_smoothness}.
\end{proof}

In particular, if either family has a nonsmooth fiber, then kernels giving a Fourier--Mukai partnership cannot be perfect.
Consequently, bounded pseudocoherent kernels, rather than perfect ones, are crucial in relative Fourier--Mukai theory.

It is not possible to replace smoothness with regularity in \Cref{thm:smoothness}.
Over a field, regularity can be detected by the equality $\operatorname{Perf}(X)=D^b_{\operatorname{coh}}(X)$.
However, over an imperfect field, a regular proper scheme need not be smooth.
For such a scheme, no Fourier--Mukai autoequivalence has a perfect kernel.
We give a concrete case.

\begin{example}
	\label{ex:nonperfect_kernels_are_needed}
	Consider any regular scheme $Y$ proper over a field $k$ that is not smooth
	(i.e.\ not geometrically regular).
	For instance, let $k:= \mathbb{F}_3 (t)$ be the rational function field in one variable, and let $Y$ be the projective plane curve over $k$ given by the equation $y^2 z + x^3 - t z^3 = 0$ \cite[Remark 16]{Kollar:2011}.
	By \Cref{prop:kernel_smoothness}, any autoequivalence of the form $\Phi_E$ on $D^b_{\operatorname{coh}}(Y)$ with $E\in D^b_{\operatorname{coh}}(Y\times_kY)$ satisfies $E\not\in \operatorname{Perf}(Y\times_kY)$.
	For example, the pushforward of $\mathcal{O}_Y$ along the diagonal, i.e.\ the kernel of the identity functor, is not perfect.
\end{example}


Next, we show that Gorensteinness is preserved under Fourier--Mukai partnerships.

\begin{theorem}
	\label{thm:gorensteinness}
	Let $f_i\colon Y_i\to S$ be proper flat morphisms of Noetherian algebraic spaces.
	Suppose $Y_1$ and $Y_2$ are Fourier--Mukai $S$-partners given by a kernel $E\in D^b_{\operatorname{coh}}(Y_1\times_S Y_2)$.
	Then $f_1$ is Gorenstein if, and only if, $f_2$ is Gorenstein.
\end{theorem}

For projective schemes over an algebraically closed field, \cite[Theorem 4.4]{Ruiperez/Hernandez/Martin/SanchodeSalas:2009} shows that Gorensteinness passes to Fourier--Mukai partners under Cohen--Macaulay and equidimensionality constraints.
Over a field, a proper variety is Gorenstein if, and only if, its category of perfect complexes admits a Serre functor.
See \cite[Lemma 6.6]{Ballard:2011} for projective varieties and \cite[Proposition 2.1]{Spence:2023} for proper varieties (i.e.\ integrality included).

\Cref{thm:gorensteinness} strengthens these statements for Gorensteinness to proper flat families of Noetherian algebraic spaces.
We do not require the $f_i$ to be projective, nor do we impose Cohen--Macaulay, equidimensionality, or integrality hypotheses.
Moreover, the result applies over any Noetherian base, including in mixed characteristic.

\subsection{Methodology}
\label{sec:intro_mehods}

For the proof of \Cref{thm:smoothness}, we base change the partnership along morphisms from fields to $S$.
This uses \cite[Theorem 10.13 and Corollary 6.18]{GuisadoVillaalgordo/Lank:2026}.
This reduces the question to regularity, which is characterized by the equality $\operatorname{Perf}=D^b_{\operatorname{coh}}$.
See \Cref{lem:regular_iff_perfect_is_dbcoh}.
For the perfectness of kernels, we use that perfectness is detected on fibers. See \Cref{lem:perf_by_fibers_of_morphism}.
Moreover, over a regular base, we show that a perfect kernel forces the partners to be regular.
See \Cref{lem:obstruction_to_singular_equivalence}.

Our proof of \Cref{thm:gorensteinness} does not pass to fibers, and it is independent of \cite{Ruiperez/Hernandez/Martin/SanchodeSalas:2009} (e.g.\ l.c.i.\ cycles are not used).
Instead, we work over the base with relative Serre functors. 
This is motivated by \cite[Definition 2.5]{SanchodeSalas/SanchodeSalas:2012}.
Specifically, we prove that a proper flat morphism is Gorenstein if, and only if, its category of perfect complexes admits a relative Serre functor.
See \Cref{lem:Gorenstein_iff_serre_functor_relative}.
This is a relative version of \cite[Lemma 6.6]{Ballard:2011}.

The relationship with diagonal morphisms with smoothness is classical: for a finite type morphism of locally Noetherian schemes, smoothness is equivalent to flatness and perfectness of the diagonal \cite[\href{https://stacks.math.columbia.edu/tag/0FDP}{Tag 0FDP}]{StacksProject}.
Lunts--Schn\"{u}rer have given a dg-categorical viewpoint which connects geometric smoothness to homological smoothness of enhancements of $\operatorname{Perf}$ \cite{Lunts/Schnurer:2016}. 
Calabrese showed a reconstruction result for algebraic spaces with affine diagonals and proper flat morphisms with Gorenstein fibers with suitable ampleness constraints \cite{Calabrese:2018}.
Our results complement these perspectives by addressing invariance of smoothness and Gorensteinness without projectivity- or ampleness-type constraints.

\begin{ack}
	The authors thank Kabeer Manali Rahul and Nebojsa Pavic for discussions.
	Lank was supported by the ERC Advanced Grant 101095900-TriCatApp.
	Guisado Villalgordo was supported under the grant PRE2022-104796 by the Spanish Ministry of Science and Innovation.
\end{ack}

\paragraph*{\textbf{AI disclosure}}

The authors used large language models (LLMs) for assistance with proofreading (e.g.\ grammar, notation, and exposition) and checking arguments (e.g.\ identifying possible gaps or errors). 
The main results and methods are due to the authors.

\section{Preliminaries}
\label{sec:prelim}

\subsection{Generation}
\label{sec:prelim_generation}

We discuss a form of generation for triangulated categories.
See \cite{Bondal/VandenBergh:2003} for details.
Let $\mathcal{T}$ be a triangulated category with shift functor $[1]\colon \mathcal{T} \to \mathcal{T}$.
Consider a subcategory $\mathcal{S} \subseteq \mathcal{T}$.
A triangulated subcategory of $\mathcal{T}$ is called \textbf{thick} if it is closed under direct summands.
Denote by $\langle \mathcal{S} \rangle$ the smallest thick subcategory of $\mathcal{T}$ containing $\mathcal{S}$.
If $\mathcal{S}$ consists of a single object $G$, we write $\langle G \rangle := \langle \mathcal{S} \rangle$.
Set $\operatorname{add}(\mathcal{S})$ to be the smallest strictly full subcategory of $\mathcal{T}$ containing $\mathcal{S}$ that is closed under shifts, finite coproducts, and direct summands.
Inductively, let $\langle \mathcal{S} \rangle_0$ consist of all zero objects in $\mathcal{T}$, $\langle \mathcal{S} \rangle_1 := \operatorname{add}(\mathcal{S})$, and
\begin{displaymath}
	\langle \mathcal{S} \rangle_n := \operatorname{add} \{ \operatorname{cone}(\phi) : \phi \in \operatorname{Hom}_{\mathcal{T}} (\langle \mathcal{S} \rangle_{n-1}, \langle \mathcal{S} \rangle_1) \}.
\end{displaymath}
It can be checked that $\langle \mathcal{S} \rangle = \bigcup_{n=0}^\infty \langle \mathcal{S} \rangle_n$.
We say $E\in \mathcal{T}$ is \textbf{finitely built by $\mathcal{S}$} if $E\in \langle \mathcal{S}\rangle$.

\begin{example}
	\label{ex:regular_ring}
	Let $X=\operatorname{Spec}(R)$ where $R$ is a regular ring.
	Then $\langle \mathcal{O}_X \rangle = D^b_{\operatorname{coh}}(X)$.
	See \cite[\href{https://stacks.math.columbia.edu/tag/066Z}{Tag 066Z}]{StacksProject} and \cite[Lemma 1.2]{Neeman:1992}.
\end{example}

\subsection{Algebraic spaces}
\label{sec:prelim_spaces}

We follow \cite{StacksProject} for conventions.
We call a morphism \textbf{natural} if it is functorial, and \textbf{canonical} if it is given by an explicit construction that need not be functorial.
Let $Y$ be a quasi-compact algebraic space.
An \textbf{\'{e}tale presentation} of $Y$ is a finitely presented, \'{e}tale, surjective morphism $U \to Y$ from a scheme.

\subsubsection{Notation}
\label{sec:prelim_notation}

We work with the small \'{e}tale ringed site $(Y_{\textrm{et}},\mathcal{O}_Y)$.
Write $D(Y)$ for the unbounded derived category of sheaves of $\mathcal{O}_Y$-modules.
Define $D_{\operatorname{qc}}(Y)$ to be the strictly full subcategory of $D(Y)$ consisting of complexes with quasi-coherent cohomology.
For $?\in \{b,+,-\}$, let $D^{?}_{\operatorname{qc}}(Y)$ be the strictly full subcategory of $D_{\operatorname{qc}}(Y)$ with prescribed cohomological vanishing.
Denote by $\operatorname{Perf}(Y)$ the strictly full subcategory of $D_{\operatorname{qc}}(Y)$ consisting of perfect complexes.
If $Y$ is Noetherian, $D^b_{\operatorname{coh}}(Y)$ is the strictly full subcategory of $D_{\operatorname{qc}}(Y)$ consisting of bounded pseudocoherent complexes.
Recall that the derived internal Hom $\operatorname{\mathbb{R}\mathcal{H}\! \mathit{om}}_{\mathcal{O}_Y}$ is formed in the ambient category $D(Y)$ \cite[\href{https://stacks.math.columbia.edu/tag/08J7}{Tag 08J7}]{StacksProject}.
In particular, it need not agree with the right adjoint to tensor product in $D_{\operatorname{qc}}(Y)$, which is denoted by $\operatorname{\mathbf{R}\mathcal{H}\! \mathit{om}}$.
See \cite[Lemma 2.11]{GuisadoVillaalgordo/Lank:2026} and \cite[\href{https://stacks.math.columbia.edu/tag/0A8A}{Tag 0A8A}]{StacksProject}. 

\subsubsection{Duality}
\label{sec:prelim_duality}

We use the duality conventions of \cite[\S 2.4]{GuisadoVillaalgordo/Lank:2026} and \cite{Neeman:2023}. 
Let $f\colon Y\to X$ be a proper morphism of Noetherian algebraic spaces.
As $f$ is proper, $f^!$ is naturally isomorphic to the right adjoint $f^\times$ of $\mathbf{R}f_{\ast}$ on quasi-coherent derived categories.
See \cite[Lemma 2.22]{GuisadoVillaalgordo/Lank:2026}. 
Recall that $f^!$ is the upper shriek functor.
It is constructed from Nagata compactifications, and defined on $D_{\operatorname{qc}}$. 
See \cite[\S 2.4]{GuisadoVillaalgordo/Lank:2026} for details.
In particular, loc.\ cit.\ discusses the results of \cite[Theorem 1.8]{Neeman:2023} for algebraic spaces, where Neeman's result constructs the so-called upper shriek functor (e.g.\ $f^!$) on $D_{\operatorname{qc}}$ and not just on $D^+_{\operatorname{qc}}$.

\subsubsection{Integral transforms}
\label{sec:prelim_integral_transforms}

We refer to \cite{GuisadoVillaalgordo/Lank:2026}.
Let $S$ be a Noetherian algebraic space. 
Suppose $Y_1$ and $Y_2$ are proper and flat algebraic $S$-spaces. 
Denote by $p_i\colon Y_1\times_S Y_2\to Y_i$ the canonical projections.
For $E\in D_{\operatorname{qc}}(Y_1\times_S Y_2)$, the \textbf{integral transform with kernel $E$} is the functor
\begin{displaymath}
    \Phi_E:=\mathbf{R}(p_2)_{\ast}(\mathbf{L}p_1^{\ast}(-)\otimes^{\mathbf{L}}E)\colon
    D_{\operatorname{qc}}(Y_1)\to D_{\operatorname{qc}}(Y_2).
\end{displaymath}
We say that $E\in D^b_{\operatorname{coh}}(Y_1\times_S Y_2)$ induces a \textbf{Fourier--Mukai equivalence} (over $S$) if $\Phi_E$ restricts to an equivalence $D^b_{\operatorname{coh}}(Y_1)\to D^b_{\operatorname{coh}}(Y_2)$.
In this case, $Y_1$ and $Y_2$ are called \textbf{Fourier--Mukai $S$-partners}, and we say it is \textbf{given} by $E$.
If clear from context, we omit the `$S$-' and say `Fourier--Mukai partners'.
We say that $A\in D_{\operatorname{qc}}(Y_1\times_S Y_2)$ is \textbf{relatively perfect over $Y_i$} if $A\otimes^{\mathbf{L}} \mathbf{L}p_i^\ast D^b_{\operatorname{qc}}(Y_i)\subseteq D^b_{\operatorname{qc}}(Y_1\times_S Y_2)$. 
By \cite[Corollary 1.4]{GuisadoVillaalgordo/Lank:2026}, a pseudocoherent complex $E$ on $Y_1\times_S Y_2$ is relatively perfect over $Y_1$ (resp.\ $Y_2$) if, and only if, $\Phi_E (D^b_{\operatorname{coh}}(Y_1))\subseteq D^b_{\operatorname{coh}}(Y_2)$ (resp.\ $\Phi_E (\operatorname{Perf}(Y_1))\subseteq \operatorname{Perf}(Y_2)$).
In both cases, any such $E$ is bounded.

\section{Auxiliary lemmas}
\label{sec:auxiliary}

\begin{lemma}
	\label{lem:invertible_via_etale_presentation}
	Let $X$ be a quasi-compact quasi-separated algebraic space.
	Let $E\in D(X)$.
	Then $E$ is invertible if, and only if, there exists an \'{e}tale presentation $s\colon U \to X$ such that $\mathbf{L}s^\ast E$ is invertible.
\end{lemma}

\begin{proof}
	If $E$ is invertible, then so is $\mathbf{L}s^\ast E$ for any \'{e}tale presentation $s$.
	Indeed, the derived pullback is monoidal \cite[\href{https://stacks.math.columbia.edu/tag/07A4}{Tag 07A4}]{StacksProject}.
	Conversely, suppose $s\colon U \to X$ is an \'{e}tale presentation such that $\mathbf{L}s^\ast E$ is invertible.
	Fix any \'{e}tale morphism $t\colon V \to X$ from a scheme.
	Consider the fibered square
	\begin{displaymath}
		\begin{tikzcd}
			{V\times_X U} & U \\
			V & {X.}
			\arrow["{t^\prime}", from=1-1, to=1-2]
			\arrow["{s^\prime}"', from=1-1, to=2-1]
			\arrow["s", from=1-2, to=2-2]
			\arrow["t"', from=2-1, to=2-2]
		\end{tikzcd}
	\end{displaymath}
	Suppose $h\colon W \to V\times_X U$ is an \'{e}tale surjective morphism from a scheme.
	Note that $s^\prime \circ h$ is an \'{e}tale covering of $V$ in $X_{\textrm{et}}$.
	As the derived pullback is monoidal, it follows that $\mathbf{L}(s \circ t^\prime \circ h)^\ast E$ is invertible.
	Recall that an algebraic space endowed with the structure on its small \'{e}tale site is a locally ringed site \cite[\href{https://stacks.math.columbia.edu/tag/04KH}{Tag 04KH}]{StacksProject}.
	Furthermore, the derived pullback along an \'{e}tale morphism from a scheme to $X$ can be identified with restriction on the small \'{e}tale site \cite[\href{https://stacks.math.columbia.edu/tag/04M4}{Tag 04M4}]{StacksProject}.
	By \cite[\href{https://stacks.math.columbia.edu/tag/0FPY}{Tag 0FPY}]{StacksProject}, there exists an \'{e}tale covering $\{g_i\colon W_i \to W\}$ by schemes such that each $E|_{W_i}$ can be represented by a shift of an invertible $\mathcal{O}_{W_i}$-module.
	Consequently, we have a covering $\{s^\prime \circ h \circ g_i\colon W_i \to V\}$ for which \cite[\href{https://stacks.math.columbia.edu/tag/0FPY}{Tag 0FPY}]{StacksProject} holds.
	Since $V$ was arbitrary, this completes the proof.
\end{proof}

\begin{lemma}
	\label{lem:upper_shriek_equals_pullback_for_etale}
	Let $f\colon Y \to X$ be a separated \'{e}tale morphism of Noetherian algebraic spaces.
	Then $f^!$ is naturally isomorphic to $\mathbf{L}f^\ast$ on $D_{\operatorname{qc}}(X)$.
\end{lemma}

\begin{proof}
	As $f$ is flat, \cite[Lemma 2.24]{GuisadoVillaalgordo/Lank:2026} gives a natural isomorphism $\mathbf{L}f^\ast E \otimes^{\mathbf{L}} f^! \mathcal{O}_X \cong f^! E$ for all $E\in D_{\operatorname{qc}}(X)$.
	Hence, it suffices to prove that $f^! \mathcal{O}_X \cong \mathcal{O}_Y$.
	Consider the fibered square
	\begin{displaymath}
		\begin{tikzcd}
			{Y\times_X Y} & Y \\
			Y & {X.}
			\arrow["{r_2}", from=1-1, to=1-2]
			\arrow["{r_1}"', from=1-1, to=2-1]
			\arrow["f", from=1-2, to=2-2]
			\arrow["f"', from=2-1, to=2-2]
		\end{tikzcd}
	\end{displaymath}
	Denote by $\Delta_f\colon Y \to Y\times_X Y$ the diagonal.
	Note that $r_1\circ \Delta_f$ and $r_2\circ \Delta_f$ are the identity of $Y$.
	By \cite[\href{https://stacks.math.columbia.edu/tag/05W1}{Tag 05W1}]{StacksProject}, $\Delta_f$ is an open immersion.
	It also is a closed immersion because $f$ is separated.
	The open immersion $\Delta_f$ admits the Nagata compactification given by
	\begin{displaymath}
		Y \xrightarrow{\Delta_f} Y\times_X Y \xrightarrow{1_{Y\times_X Y}}Y\times_X Y,
	\end{displaymath}
	it follows that $\mathbf{L}\Delta_f^\ast\cong \Delta_f^!$.
	There exists a string of isomorphisms
	\begin{displaymath}
		\begin{aligned}
			f^! \mathcal{O}_X
			&\cong \mathbf{L}\Delta_f^\ast \mathbf{L}r_2^\ast f^! \mathcal{O}_X && (r_2\circ \Delta_f = 1_Y)
			\\&\cong \mathbf{L}\Delta_f^\ast r_1^! \mathcal{O}_Y && (\textrm{\cite[Lemma 2.23]{GuisadoVillaalgordo/Lank:2026}})
			\\&\cong \Delta_f^! r_1^! \mathcal{O}_Y
			\\&\cong (r_1 \circ \Delta_f)^! \mathcal{O}_Y && (\textrm{\cite[Lemma 2.21]{GuisadoVillaalgordo/Lank:2026}})
			\\&\cong \mathcal{O}_Y && (r_1\circ \Delta_f = 1_Y).
		\end{aligned}
	\end{displaymath}
	This completes the proof.
\end{proof}

\begin{lemma}
	\label{lem:pullback_dualizing_complex}
	Let $f\colon Y \to X$ be a separated morphism of finite type between Noetherian algebraic spaces.
	If $K$ is a dualizing complex on $X$, then $f^! K$ is a dualizing complex on $Y$.
\end{lemma}

\begin{proof}
	Choose \'{e}tale presentations $t\colon V \to Y$ and $s\colon U \to X$ from affine schemes.
	There exists a commutative diagram
	\begin{displaymath}
		\begin{tikzcd}
			{V\times_X U} & {Y\times_X U} & U \\
			V & Y & X
			\arrow["{t^\prime}", from=1-1, to=1-2]
			\arrow["{s^{\prime\prime}}"', from=1-1, to=2-1]
			\arrow["{f^\prime}", from=1-2, to=1-3]
			\arrow["{s^\prime}", from=1-2, to=2-2]
			\arrow["s", from=1-3, to=2-3]
			\arrow["t"', from=2-1, to=2-2]
			\arrow["f"', from=2-2, to=2-3]
		\end{tikzcd}
	\end{displaymath}
	obtained by base change.
	By \cite[Lemma 2.19]{GuisadoVillaalgordo/Lank:2026}, $t$, $s$, and $f\circ t$ are quasi-affine.
	Base change implies $s^\prime$, $s^{\prime\prime}$, and $t^\prime$ are quasi-affine, finitely presented, and \'{e}tale.
	In particular, $V\times_X U$ is a Noetherian scheme.
	There is a direct computation,
	\begin{displaymath}
		\begin{aligned}
			\mathbf{L}(s^{\prime\prime})^\ast \mathbf{L}t^\ast f^! K
			&\cong \mathbf{L}(s^{\prime\prime})^\ast t^! f^! K && (\textrm{\Cref{lem:upper_shriek_equals_pullback_for_etale}})
			\\&\cong \mathbf{L}(s^{\prime\prime})^\ast (f\circ t)^! K && (\textrm{\cite[Lemma 2.21]{GuisadoVillaalgordo/Lank:2026}})
			\\&\cong (f^\prime \circ t^\prime)^! \mathbf{L}s^\ast K && (\textrm{\cite[Lemma 2.23]{GuisadoVillaalgordo/Lank:2026}}).
		\end{aligned}
	\end{displaymath}
	Note that $\mathbf{L}s^\ast K$ is a dualizing complex.
	As $f^\prime \circ t^\prime$ is a separated finite type morphism of Noetherian schemes, \cite[\href{https://stacks.math.columbia.edu/tag/0AA3}{Tag 0AA3}]{StacksProject} says that $(f^\prime \circ t^\prime)^! \mathbf{L}s^\ast K$ is a dualizing complex.
	Since $t\circ s^{\prime\prime}$ is an \'{e}tale presentation of $Y$,  $f^! K$ is a dualizing complex \cite[\href{https://stacks.math.columbia.edu/tag/0E4Y}{Tag 0E4Y}]{StacksProject}.
\end{proof}

\section{Gorensteinness}
\label{sec:Gorensteinness}

\begin{definition}
	Let $f\colon Y \to X$ be a proper flat morphism of Noetherian algebraic spaces.
	An \textbf{$f$-Serre functor} is an autoequivalence $\mathbb{S}_f\colon \operatorname{Perf}(Y) \to \operatorname{Perf}(Y)$ equipped with a bifunctorial isomorphism
	\begin{displaymath}
		\operatorname{\mathbf{R}\mathcal{H}\! \mathit{om}} (\mathbf{R}f_\ast \operatorname{\mathbf{R}\mathcal{H}\! \mathit{om}} (P,Q),\mathcal{O}_X) \cong \mathbf{R}f_\ast \operatorname{\mathbf{R}\mathcal{H}\! \mathit{om}} (Q,\mathbb{S}_f (P))
	\end{displaymath}
	for all $P,Q\in \operatorname{Perf}(Y)$.
\end{definition}

\begin{remark}
	The definition above is an analog of \cite[Definition 2.5]{SanchodeSalas/SanchodeSalas:2012}.
\end{remark}

\begin{lemma}
	\label{lem:Gorenstein_morphism_iff_invertible}
	Let $f\colon Y \to X$ be a proper flat morphism of Noetherian algebraic spaces.
	Then $f$ is Gorenstein if, and only if, $f^!\mathcal{O}_X$ is invertible.
\end{lemma}

\begin{proof}
	Choose an \'{e}tale presentation $s\colon U \to X$ from an affine scheme.
	Note that $U$ is Noetherian.
	Consider the fibered square
	\begin{displaymath}
		\begin{tikzcd}
			{Y\times_X U} & U \\
			Y & {X.}
			\arrow["{f^\prime}", from=1-1, to=1-2]
			\arrow["{s^\prime}"', from=1-1, to=2-1]
			\arrow["s", from=1-2, to=2-2]
			\arrow["f"', from=2-1, to=2-2]
		\end{tikzcd}
	\end{displaymath}
	Choose an \'{e}tale presentation $w\colon V \to Y\times_X U$ from an affine scheme.
	Note that $V$ is Noetherian, and that $w$ is separated by \cite[Lemma 2.19]{GuisadoVillaalgordo/Lank:2026}.
	The composition $f^\prime \circ w\colon V \to U$ is a flat separated finite type morphism of Noetherian schemes.
	There is a string of isomorphisms,
	\begin{displaymath}
		\begin{aligned}
			(f^\prime \circ w)^! \mathcal{O}_U
			&\cong w^! (f^\prime)^! \mathcal{O}_U && (\textrm{\cite[Lemma 2.21]{GuisadoVillaalgordo/Lank:2026}})
			\\&\cong w^! \mathbf{L}(s^\prime)^\ast f^! \mathcal{O}_X && (\textrm{\cite[Lemma 2.23]{GuisadoVillaalgordo/Lank:2026}})
			\\&\cong \mathbf{L}w^\ast \mathbf{L}(s^\prime)^\ast f^! \mathcal{O}_X && (\textrm{\Cref{lem:upper_shriek_equals_pullback_for_etale}})
			\\&\cong \mathbf{L}(s^\prime \circ w)^\ast f^! \mathcal{O}_X.
		\end{aligned}
	\end{displaymath}
	Since $s^\prime \circ w$ is an \'{e}tale presentation, \Cref{lem:invertible_via_etale_presentation} implies that $f^! \mathcal{O}_X$ is invertible if, and only if, $(f^\prime \circ w)^! \mathcal{O}_U$ is invertible.
	By \cite[\href{https://stacks.math.columbia.edu/tag/0C08}{Tag 0C08}]{StacksProject}, $f^\prime \circ w$ is Gorenstein if, and only if, $(f^\prime \circ w)^! \mathcal{O}_U$ is isomorphic Zariski locally on $V$ to an invertible object.
	By \cite[\href{https://stacks.math.columbia.edu/tag/0FPY}{Tag 0FPY}]{StacksProject}, the latter holds if, and only if, $(f^\prime \circ w)^! \mathcal{O}_U$ is invertible.
	Now, $s$ and $s^\prime \circ w$ are \'{e}tale, and $s^\prime \circ w$ is surjective.
	Hence, by \cite[\href{https://stacks.math.columbia.edu/tag/0E14}{Tags 0E14} and \href{https://stacks.math.columbia.edu/tag/04NC}{04NC}]{StacksProject}, $f$ is Gorenstein if, and only if, $f^\prime \circ w$ is Gorenstein.
	This completes the proof.
\end{proof}

\begin{lemma}
	\label{lem:Gorenstein_iff_serre_functor_relative}
	Let $f\colon Y \to X$ be a proper flat morphism of Noetherian algebraic spaces.
	Then $\operatorname{Perf}(Y)$ admits an $f$-Serre functor if, and only if, $f$ is Gorenstein.
\end{lemma}

\begin{proof}
	There is a string of bifunctorial isomorphisms for all $P,Q\in \operatorname{Perf}(Y)$:
	\begin{displaymath}
		\begin{aligned}
			& \operatorname{\mathbf{R}\mathcal{H}\! \mathit{om}}
			(\mathbf{R}f_\ast \operatorname{\mathbf{R}\mathcal{H}\! \mathit{om}} (P,Q),\mathcal{O}_X)
			\\&\cong \mathbf{R}f_\ast  \operatorname{\mathbf{R}\mathcal{H}\! \mathit{om}} (\operatorname{\mathbf{R}\mathcal{H}\! \mathit{om}} (P,Q), f^\times \mathcal{O}_X) && (\textrm{\cite[Lemma 5.9]{GuisadoVillaalgordo/Lank:2026}})
			\\&\cong \mathbf{R}f_\ast  \operatorname{\mathbf{R}\mathcal{H}\! \mathit{om}} (\operatorname{\mathbf{R}\mathcal{H}\! \mathit{om}} (P,Q), f^! \mathcal{O}_X) && (\textrm{\cite[Lemma 2.22]{GuisadoVillaalgordo/Lank:2026}})
			\\&\cong \mathbf{R}f_\ast  \operatorname{\mathbf{R}\mathcal{H}\! \mathit{om}} (\operatorname{\mathbf{R}\mathcal{H}\! \mathit{om}} (P,\mathcal{O}_Y) \otimes^{\mathbf{L}} Q , f^! \mathcal{O}_X) && (\textrm{\cite[\href{https://stacks.math.columbia.edu/tag/08JJ}{Tag 08JJ} and \href{https://stacks.math.columbia.edu/tag/0D0S}{0D0S}]{StacksProject}})
			\\&\cong \mathbf{R}f_\ast  \operatorname{\mathbf{R}\mathcal{H}\! \mathit{om}} ( Q , \operatorname{\mathbf{R}\mathcal{H}\! \mathit{om}} (\operatorname{\mathbf{R}\mathcal{H}\! \mathit{om}} (P,\mathcal{O}_Y) , f^! \mathcal{O}_X ) ) && (\textrm{\cite[\href{https://stacks.math.columbia.edu/tag/08J9}{Tag 08J9} and \href{https://stacks.math.columbia.edu/tag/0D0S}{0D0S}]{StacksProject}})
			\\&\cong \mathbf{R}f_\ast  \operatorname{\mathbf{R}\mathcal{H}\! \mathit{om}} ( Q , \operatorname{\mathbf{R}\mathcal{H}\! \mathit{om}} (\operatorname{\mathbf{R}\mathcal{H}\! \mathit{om}} (P,\mathcal{O}_Y) , \mathcal{O}_Y ) \otimes^{\mathbf{L}} f^! \mathcal{O}_X ) && (\textrm{\cite[\href{https://stacks.math.columbia.edu/tag/08JJ}{Tag 08JJ} and \href{https://stacks.math.columbia.edu/tag/0D0S}{0D0S}]{StacksProject}})
			\\&\cong \mathbf{R}f_\ast  \operatorname{\mathbf{R}\mathcal{H}\! \mathit{om}} ( Q , P \otimes^{\mathbf{L}} f^! \mathcal{O}_X ).
		\end{aligned}
	\end{displaymath}
	If $f$ is Gorenstein, then \Cref{lem:Gorenstein_morphism_iff_invertible} implies $f^! \mathcal{O}_X$ is invertible.
	In such a case, $(-)\otimes^{\mathbf{L}} f^! \mathcal{O}_X\colon \operatorname{Perf}(Y) \to \operatorname{Perf}(Y)$ is an autoequivalence, and the bifunctorial isomorphisms above show that this defines an $f$-Serre functor.

	Conversely, if $\mathbb{S}_f$ is an $f$-Serre functor on $\operatorname{Perf}(Y)$, then for all $P,Q\in \operatorname{Perf}(Y)$ there is a string of natural isomorphisms:
	\begin{displaymath}
		\begin{aligned}
			\operatorname{Hom}(Q,\mathbb{S}_f (P))
			&\cong \operatorname{Hom}(\mathcal{O}_Y , \operatorname{\mathbf{R}\mathcal{H}\! \mathit{om}}(Q,\mathbb{S}_f (P))) && (\textrm{\cite[\href{https://stacks.math.columbia.edu/tag/08J7}{Tag 08J7} and \href{https://stacks.math.columbia.edu/tag/0D0S}{0D0S}]{StacksProject}})
			\\&\cong \operatorname{Hom}(\mathbf{L}f^\ast \mathcal{O}_X , \operatorname{\mathbf{R}\mathcal{H}\! \mathit{om}}(Q,\mathbb{S}_f (P)))
			\\&\cong \operatorname{Hom}(\mathcal{O}_X , \mathbf{R}f_\ast \operatorname{\mathbf{R}\mathcal{H}\! \mathit{om}}(Q,\mathbb{S}_f (P))) && \textrm{(pull/push)}
			\\&\cong \operatorname{Hom}(\mathcal{O}_X , \operatorname{\mathbf{R}\mathcal{H}\! \mathit{om}} (\mathbf{R}f_\ast \operatorname{\mathbf{R}\mathcal{H}\! \mathit{om}} (P,Q),\mathcal{O}_X) ) && \textrm{(Serre functor)}
			\\&\cong \operatorname{Hom}(\mathbf{R}f_\ast \operatorname{\mathbf{R}\mathcal{H}\! \mathit{om}}(P,Q), \mathcal{O}_X )
			\\&\cong \operatorname{Hom}(\operatorname{\mathbf{R}\mathcal{H}\! \mathit{om}}(P,Q), f^\times \mathcal{O}_X ) && (\textrm{right adjoint of push})
			\\&\cong \operatorname{Hom}(\operatorname{\mathbf{R}\mathcal{H}\! \mathit{om}}(P,Q), f^! \mathcal{O}_X ) && \textrm{(\cite[Lemma 2.22]{GuisadoVillaalgordo/Lank:2026})}
			\\&\cong \operatorname{Hom}(Q,P\otimes^{\mathbf{L}} f^! \mathcal{O}_X).
		\end{aligned}
	\end{displaymath}
	Fix $P\in \operatorname{Perf}(Y)$.
	As $\mathbb{S}_f(P)$ is perfect, \cite[\href{https://stacks.math.columbia.edu/tag/001P}{Tag 001P}]{StacksProject}, applied to these isomorphisms as functors of $Q\in\operatorname{Perf}(Y)$, yields
	\begin{displaymath}
		\phi_P\colon \mathbb{S}_f (P) \to P\otimes^{\mathbf{L}} f^! \mathcal{O}_X
	\end{displaymath}
	such that $\operatorname{Hom}(Q,\phi_P)$ is an isomorphism for all $Q\in \operatorname{Perf}(Y)$.
	Pick a compact generator $G$ for $D_{\operatorname{qc}}(Y)$, see \cite[\href{https://stacks.math.columbia.edu/tag/09IY}{Tag 09IY}]{StacksProject}.
	Apply $\operatorname{Hom}(G,-)$ to the distinguished triangle
	\begin{displaymath}
		\mathbb{S}_f (P) \xrightarrow{\phi_P} P \otimes^{\mathbf{L}} f^! \mathcal{O}_X \to \operatorname{cone}(\phi_P) \to \mathbb{S}_f (P)[1].
	\end{displaymath}
	There exists a long exact sequence (see \cite[\href{https://stacks.math.columbia.edu/tag/0149}{Tag 0149}]{StacksProject}),
	\begin{displaymath}
		\cdots \to \operatorname{Hom}(G,\mathbb{S}_f (P)) \xrightarrow{\operatorname{Hom}(G,\phi_P)} \operatorname{Hom}(G,P \otimes^{\mathbf{L}} f^! \mathcal{O}_X) \to \operatorname{Hom}(G,\operatorname{cone}(\phi_P)) \to \cdots.
	\end{displaymath}
	As $\operatorname{Hom}(G[n],\phi_P)$ is an isomorphism for all $n\in\mathbb{Z}$, it follows that $\operatorname{Hom}(G,\operatorname{cone}(\phi_P)[n])\cong 0$ for all $n\in \mathbb{Z}$.
	Since $G$ compactly generates $D_{\operatorname{qc}}(Y)$, we obtain that $\operatorname{cone}(\phi_P)\cong 0$.
	In other words, $\phi_P$ is an isomorphism.
	Taking $P:=\mathcal{O}_Y$ implies $f^! \mathcal{O}_X$ is perfect.
	Since $\mathbb{S}_f$ is essentially surjective, there exists $P_0 \in \operatorname{Perf}(Y)$ such that $\mathbb{S}_f (P_0) \cong \mathcal{O}_Y$.
	Then
	\begin{displaymath}
		\mathcal{O}_Y \cong \mathbb{S}_f (P_0) \cong P_0 \otimes^{\mathbf{L}} f^! \mathcal{O}_X.
	\end{displaymath}
	Since $P_0$ and $f^! \mathcal{O}_X$ are perfect, the preceding isomorphism shows that $f^! \mathcal{O}_X$ is invertible, with inverse $P_0$.
	Therefore, by \Cref{lem:Gorenstein_morphism_iff_invertible}, $f$ is Gorenstein.
\end{proof}

\begin{lemma}
	\label{lem:Gorenstein_alt_proof1}
	Let $f_i\colon Y_i \to S$ be proper flat morphisms of Noetherian algebraic spaces.
	Suppose $Y_1$ and $Y_2$ are Fourier--Mukai $S$-partners given by $B\in D^b_{\operatorname{coh}}(Y_1\times_S Y_2)$.
	There are natural isomorphisms for all $P,Q\in \operatorname{Perf}(Y_1)$:
	\begin{enumerate}
		\item $\mathbf{R}(f_1)_\ast \operatorname{\mathbf{R}\mathcal{H}\! \mathit{om}} (P,Q)
		\cong
		\mathbf{R}(f_2)_\ast \operatorname{\mathbf{R}\mathcal{H}\! \mathit{om}} (\Phi_B (P), \Phi_B (Q))$
		\item $\operatorname{\mathbf{R}\mathcal{H}\! \mathit{om}} ( \mathbf{R}(f_1)_\ast \operatorname{\mathbf{R}\mathcal{H}\! \mathit{om}} (P,Q) ,\mathcal{O}_S)
		\cong
		\operatorname{\mathbf{R}\mathcal{H}\! \mathit{om}} (\mathbf{R}(f_2)_\ast \operatorname{\mathbf{R}\mathcal{H}\! \mathit{om}} (\Phi_B (P), \Phi_B (Q)), \mathcal{O}_S)$.
	\end{enumerate}
\end{lemma}

\begin{proof}
	We prove the first claim.
	By \cite[Remark 2.7 and Lemma 10.9]{GuisadoVillaalgordo/Lank:2026}, $\Phi_B$ induces an equivalence on $D_{\operatorname{qc}}$.
	Denote by $p_i\colon Y_1\times_S Y_2\to Y_i$ the projections.
	For $P,Q\in\operatorname{Perf}(Y_1)$ and $E\in D_{\operatorname{qc}}(S)$, there exist a string of natural isomorphisms:
	\begin{displaymath}
		\begin{aligned}
			&\operatorname{Hom}(E, \mathbf{R}(f_1)_\ast \operatorname{\mathbf{R}\mathcal{H}\! \mathit{om}} (P,Q))
			\\&\cong \operatorname{Hom}(\mathbf{L}f_1^\ast  E, \operatorname{\mathbf{R}\mathcal{H}\! \mathit{om}} (P,Q)) && (\textrm{pull/push})
			\\&\cong \operatorname{Hom}(\mathbf{L}f_1^\ast  E \otimes^{\mathbf{L}} P,  Q) && (\textrm{\cite[\href{https://stacks.math.columbia.edu/tag/08J7}{Tag 08J7} and \href{https://stacks.math.columbia.edu/tag/0D0S}{0D0S}]{StacksProject}})
			\\&\cong \operatorname{Hom}(\Phi_B (\mathbf{L}f_1^\ast  E \otimes^{\mathbf{L}} P),  \Phi_B (Q)) && (\Phi_B \textrm{ is an equivalence})
			\\&\cong \operatorname{Hom}(\mathbf{R}(p_2)_\ast (\mathbf{L}p_1^\ast (\mathbf{L}f_1^\ast  E \otimes^{\mathbf{L}} P) \otimes^{\mathbf{L}} B),  \Phi_B (Q)) && (\textrm{def.\ of }  \Phi_B)
			\\&\cong \operatorname{Hom}(\mathbf{R}(p_2)_\ast ( \mathbf{L}p_1^\ast \mathbf{L}f_1^\ast  E \otimes^{\mathbf{L}} \mathbf{L}p_1^\ast P \otimes^{\mathbf{L}} B),  \Phi_B (Q)) && \textrm{(\cite[\href{https://stacks.math.columbia.edu/tag/07A4}{Tag 07A4}]{StacksProject})}
			\\&\cong \operatorname{Hom}(\mathbf{R}(p_2)_\ast ( \mathbf{L}p_2^\ast \mathbf{L}f_2^\ast  E \otimes^{\mathbf{L}} \mathbf{L}p_1^\ast P \otimes^{\mathbf{L}} B),  \Phi_B (Q)) && (f_1 \circ p_1 = f_2 \circ p_2)
			\\&\cong \operatorname{Hom}(\mathbf{R}(p_2)_\ast ( \mathbf{L}p_1^\ast P \otimes^{\mathbf{L}} B) \otimes^{\mathbf{L}} \mathbf{L}f_2^\ast  E ,  \Phi_B (Q)) && (\textrm{projection formula})
			\\&\cong \operatorname{Hom}(\Phi_B (P) \otimes^{\mathbf{L}} \mathbf{L}f_2^\ast  E ,  \Phi_B (Q)) && (\textrm{def.\ of }  \Phi_B)
			\\&\cong \operatorname{Hom}(\mathbf{L}f_2^\ast  E, \operatorname{\mathbf{R}\mathcal{H}\! \mathit{om}}(\Phi_B (P), \Phi_B (Q))) && (\textrm{\cite[\href{https://stacks.math.columbia.edu/tag/08J7}{Tag 08J7} and \href{https://stacks.math.columbia.edu/tag/0D0S}{0D0S}]{StacksProject}})
			\\&\cong \operatorname{Hom}(E, \mathbf{R}(f_2)_\ast \operatorname{\mathbf{R}\mathcal{H}\! \mathit{om}}(\Phi_B (P), \Phi_B (Q))) && (\textrm{pull/push}).
		\end{aligned}
	\end{displaymath}
	These isomorphisms are natural in $E\in D_{\operatorname{qc}}(S)$ and in $P,Q\in \operatorname{Perf}(Y_1)$.
	Hence, \cite[\href{https://stacks.math.columbia.edu/tag/001P}{Tag 001P}]{StacksProject} gives an isomorphism
	\begin{displaymath}
		\mathbf{R}(f_1)_\ast \operatorname{\mathbf{R}\mathcal{H}\! \mathit{om}} (P,Q)
		\to \mathbf{R}(f_2)_\ast \operatorname{\mathbf{R}\mathcal{H}\! \mathit{om}}(\Phi_B (P), \Phi_B (Q))
	\end{displaymath}
	which is functorial in $P$ and $Q$.
	This proves the first claim.
	Applying $\operatorname{\mathbf{R}\mathcal{H}\! \mathit{om}}(-,\mathcal{O}_S)$ gives the second claim, and completes the proof.
\end{proof}

\begin{lemma}
	\label{lem:Gorenstein_alt_proof2}
	Let $f_i\colon Y_i \to S$ be proper flat morphisms of Noetherian algebraic spaces.
	Suppose $Y_1$ and $Y_2$ are Fourier--Mukai $S$-partners given by $B\in D^b_{\operatorname{coh}}(Y_1\times_S Y_2)$.
	Denote by $p_i\colon Y_1\times_S Y_2 \to Y_i$ the projections.
	If $\mathbb{S}_1$ is an $f_1$-Serre functor, then there are natural isomorphisms for all $P,Q\in \operatorname{Perf}(Y_1)$:
	\begin{displaymath}
		\mathbf{R}(f_2)_\ast \operatorname{\mathbf{R}\mathcal{H}\! \mathit{om}} (\Phi_B (Q), \mathbb{S}_2 (\Phi_B (P)))
		\cong \mathbf{R}(f_1)_\ast \operatorname{\mathbf{R}\mathcal{H}\! \mathit{om}} (Q, \mathbb{S}_1 (P))
	\end{displaymath}
	where $\mathbb{S}_2 := \Phi_B \circ \mathbb{S}_1 \circ \Phi_{\operatorname{\mathbf{R}\mathcal{H}\! \mathit{om}}(B,p_2^! \mathcal{O}_{Y_2})}$.
	Here, $\Phi_{\operatorname{\mathbf{R}\mathcal{H}\! \mathit{om}}(B,p_2^! \mathcal{O}_{Y_2})}$ denotes the integral transform from $D_{\operatorname{qc}}(Y_2)$ to $D_{\operatorname{qc}}(Y_1)$, i.e.\ the roles of $p_1$ and $p_2$ are interchanged.
\end{lemma}

\begin{proof}
	Fix $P\in \operatorname{Perf}(Y_1)$.
	By \cite[Remark 2.7 and Lemma 10.9]{GuisadoVillaalgordo/Lank:2026}, $\Phi_B$ is an equivalence on $D_{\operatorname{qc}}$ which restricts to an equivalence on $\operatorname{Perf}$.
	By \cite[Proposition 9.3]{GuisadoVillaalgordo/Lank:2026}, applied with $t$ the identity, $\Phi_{\operatorname{\mathbf{R}\mathcal{H}\! \mathit{om}}(B,p_2^! \mathcal{O}_{Y_2})}$ is right adjoint to $\Phi_B$ on $D_{\operatorname{qc}}$.
	Hence, it is a quasi-inverse of $\Phi_B$.
	We claim that $\mathbb{S}_2 (\Phi_B (P))\cong \Phi_B (\mathbb{S}_1 (P))$.
	Indeed,
	\begin{displaymath}
		\begin{aligned}
			\mathbb{S}_2 (\Phi_B (P))
			&\cong (\Phi_B \circ \mathbb{S}_1 \circ \Phi_{\operatorname{\mathbf{R}\mathcal{H}\! \mathit{om}}(B,p_2^! \mathcal{O}_{Y_2})} \circ \Phi_B) (P) && (\textrm{def.\ of } \mathbb{S}_2)
			\\& \cong (\Phi_B \circ \mathbb{S}_1) (P) && (\textrm{quasi-inverse}).
		\end{aligned}
	\end{displaymath}
	Then \Cref{lem:Gorenstein_alt_proof1} applied to $Q$ and $\mathbb{S}_1 (P)$ yields
	\begin{displaymath}
		\mathbf{R}(f_1)_\ast \operatorname{\mathbf{R}\mathcal{H}\! \mathit{om}} (Q,\mathbb{S}_1 (P))
		\cong
		\mathbf{R}(f_2)_\ast \operatorname{\mathbf{R}\mathcal{H}\! \mathit{om}} (\Phi_B (Q), \mathbb{S}_2 (\Phi_B (P))).
	\end{displaymath}
	Thus, the desired claim follows.
\end{proof}

\begin{proof}
	[Proof of \Cref{thm:gorensteinness}]
	\label{proof:prop:Gorenstein_derived_invariance2}
	By symmetry, it suffices to show that if $f_1$ is Gorenstein, then so is $f_2$.
	By \cite[Remark 2.7]{GuisadoVillaalgordo/Lank:2026}, $E$ is relatively perfect over each $Y_i$.
	As $f_1$ is Gorenstein, \Cref{lem:Gorenstein_iff_serre_functor_relative} gives an $f_1$-Serre functor $\mathbb{S}_1$.
	Define $\mathbb{S}_2 := \Phi_E \circ \mathbb{S}_1 \circ \Phi_{\operatorname{\mathbf{R}\mathcal{H}\! \mathit{om}}(E,p_2^! \mathcal{O}_{Y_2})}$ as in \Cref{lem:Gorenstein_alt_proof2}.
	By \cite[Lemma 10.9]{GuisadoVillaalgordo/Lank:2026}, $\Phi_E$ restricts to an equivalence $\operatorname{Perf}(Y_1)\to\operatorname{Perf}(Y_2)$, and so $\mathbb{S}_2$ is an autoequivalence of $\operatorname{Perf}(Y_2)$.
	For all $P,Q\in \operatorname{Perf}(Y_1)$, there is a string of natural isomorphisms
	\begin{displaymath}
		\begin{aligned}
			&\operatorname{\mathbf{R}\mathcal{H}\! \mathit{om}} (\mathbf{R}(f_2)_\ast \operatorname{\mathbf{R}\mathcal{H}\! \mathit{om}} (\Phi_E (P), \Phi_E (Q)), \mathcal{O}_S)
			\\&\cong \operatorname{\mathbf{R}\mathcal{H}\! \mathit{om}} ( \mathbf{R}(f_1)_\ast \operatorname{\mathbf{R}\mathcal{H}\! \mathit{om}} (P,Q) ,\mathcal{O}_S) && (\textrm{\Cref{lem:Gorenstein_alt_proof1}})
			\\&\cong \mathbf{R}(f_1)_\ast \operatorname{\mathbf{R}\mathcal{H}\! \mathit{om}} (Q, \mathbb{S}_1 (P)) && (\mathbb{S}_1 \textrm{ is an } f_1\textrm{-Serre functor})
			\\&\cong \mathbf{R}(f_2)_\ast \operatorname{\mathbf{R}\mathcal{H}\! \mathit{om}} (\Phi_E (Q), \mathbb{S}_2 (\Phi_E (P))) && (\textrm{\Cref{lem:Gorenstein_alt_proof2}}).
		\end{aligned}
	\end{displaymath}
	Note that every object of $\operatorname{Perf}(Y_2)$ is isomorphic to one of the form $\Phi_E(P)$ with $P\in \operatorname{Perf}(Y_1)$.
	Hence, from the equivalence $\Phi_E$, it follows that $\mathbb{S}_2$ is an $f_2$-Serre functor.
	Consequently, \Cref{lem:Gorenstein_iff_serre_functor_relative} implies $f_2$ is Gorenstein.
\end{proof}

\section{Smoothness}
\label{sec:smoothness}

\begin{lemma}
	\label{lem:smooth_if_all_representatives}
	Let $f\colon Y \to X$ be a finitely presented flat morphism of Noetherian algebraic spaces.
	If the base change of $f$ along every $\operatorname{Spec}(k)\to X$ is smooth where $k$ is a field, then $f$ is smooth.
\end{lemma}

\begin{proof}
	Choose an \'{e}tale presentation $s\colon U \to X$.
	Consider the fibered square
	\begin{displaymath}
		\begin{tikzcd}
			{Y\times_X U} & U \\
			Y & {X.}
			\arrow["{f^\prime}", from=1-1, to=1-2]
			\arrow["{s^\prime}"', from=1-1, to=2-1]
			\arrow["s", from=1-2, to=2-2]
			\arrow["f"', from=2-1, to=2-2]
		\end{tikzcd}
	\end{displaymath}
	Let $t\colon V \to Y\times_X U$ be an \'{e}tale presentation.
	Fix $u \in U$.
	Denote by $i_u\colon \operatorname{Spec}(\kappa(u)) \to U$ the canonical morphism.
	There exists a fibered square
	\begin{displaymath}
		\begin{tikzcd}
			{Y\times_X \operatorname{Spec}(\kappa(u))} & {\operatorname{Spec}(\kappa(u))} \\
			{Y\times_X U} & {U.}
			\arrow["{f^\prime_u}", from=1-1, to=1-2]
			\arrow["{i^\prime_u}"', from=1-1, to=2-1]
			\arrow["{i_u}", from=1-2, to=2-2]
			\arrow["{f^\prime}"', from=2-1, to=2-2]
		\end{tikzcd}
	\end{displaymath}
	By the hypothesis, $f^\prime_u$ is smooth.
	Moreover, we have another fibered square,
	\begin{displaymath}
		\begin{tikzcd}
			{V\times_U \operatorname{Spec}(\kappa(u))} & {Y\times_X \operatorname{Spec}(\kappa(u))} & {\operatorname{Spec}(\kappa(u))} \\
			V & {Y\times_X U} & {U.}
			\arrow["{t_u}", from=1-1, to=1-2]
			\arrow["{i^{\prime\prime}_u}"', from=1-1, to=2-1]
			\arrow["{f^\prime_u}", from=1-2, to=1-3]
			\arrow["{i^\prime_u}"', from=1-2, to=2-2]
			\arrow["{i_u}"', from=1-3, to=2-3]
			\arrow["t"', from=2-1, to=2-2]
			\arrow["{f^\prime}"', from=2-2, to=2-3]
		\end{tikzcd}
	\end{displaymath}
	Base change says $t_u$ is an \'{e}tale presentation.
	Hence, the composition $f^\prime_u \circ t_u$ is a smooth finitely presented morphism of Noetherian schemes.
	By \cite[\href{https://stacks.math.columbia.edu/tag/01V8}{Tag 01V8}]{StacksProject}, $f^\prime \circ t$ is smooth.
	Appealing to \cite[\href{https://stacks.math.columbia.edu/tag/03ZF}{Tag 03ZF}]{StacksProject}, we have that $f$ is smooth.
\end{proof}

\begin{lemma}
	\label{lem:smooth_if_some_representatives}
	Let $f\colon Y \to X$ be a finitely presented flat morphism of Noetherian algebraic spaces.
	If for every $p\in |X|$ there exists a representative $\operatorname{Spec}(k) \to X$ of $p$ such that the base change of $f$ along $\operatorname{Spec}(k)\to X$ is smooth, then the base change of $f$ along every $\operatorname{Spec}(k)\to X$ is smooth where $k$ is a field.
\end{lemma}

\begin{proof}
	Fix $p\in |X|$.
	Choose a representative $t\colon \operatorname{Spec}(k) \to X$ of $p$ such that the base change of $f$ along $t$ is smooth.
	Suppose $t^\prime\colon \operatorname{Spec}(k^\prime) \to X$ is any representative of $p$.
	By the definition of the points of an algebraic space \cite[\href{https://stacks.math.columbia.edu/tag/03BT}{Tag 03BT}]{StacksProject}, there exist a field $\ell$ and a commutative diagram
	\begin{displaymath}
		\begin{tikzcd}
			{\operatorname{Spec}(\ell)} & {\operatorname{Spec}(k)} \\
			{\operatorname{Spec}(k^\prime)} & {X.}
			\arrow["{q_2}", from=1-1, to=1-2]
			\arrow["{q_1}"', from=1-1, to=2-1]
			\arrow["t", from=1-2, to=2-2]
			\arrow["{t^\prime}"', from=2-1, to=2-2]
		\end{tikzcd}
	\end{displaymath}
	Consider the fibered cube
	\begin{displaymath}
		\begin{tikzcd}
			{Y\times_X \operatorname{Spec}(\ell)} &&& {Y\times_X \operatorname{Spec}(k)} && \\
			& {Y\times_X \operatorname{Spec}(k^\prime)} &&&& Y \\
			{\operatorname{Spec}(\ell)} &&& {\operatorname{Spec}(k)} \\
			& {\operatorname{Spec}(k^\prime)} &&&& {X.}
			\arrow["{q^\prime_2}", from=1-1, to=1-4]
			\arrow["{q_1^\prime}", from=1-1, to=2-2]
			\arrow["{f_\ell}", from=1-1, to=3-1]
			\arrow["{g_k}", from=1-4, to=2-6]
			\arrow["{f_k}"{pos=0.8}, from=1-4, to=3-4]
			\arrow["{g_{k^\prime}}"{pos=0.2}, from=2-2, to=2-6]
			\arrow["{f_{k^\prime}}"{pos=0.8}, from=2-2, to=4-2]
			\arrow["f", from=2-6, to=4-6]
			\arrow["{q_2}"{pos=0.7}, from=3-1, to=3-4]
			\arrow["{q_1}"', from=3-1, to=4-2]
			\arrow["t", from=3-4, to=4-6]
			\arrow["{t^\prime}"', from=4-2, to=4-6]
		\end{tikzcd}
	\end{displaymath}
	As $f_k$ is smooth, base change implies $f_\ell$ is smooth.
	Choose an \'{e}tale presentation $s\colon U \to Y$.
	There exists a further fibered cube
	\begin{displaymath}
		\begin{tikzcd}
			{U\times_X \operatorname{Spec}(\ell)} &&& {U\times_X \operatorname{Spec}(k)} && \\
			& {U\times_X \operatorname{Spec}(k^\prime)} &&&& U \\
			{Y\times_X \operatorname{Spec}(\ell)} &&& {Y\times_X \operatorname{Spec}(k)} \\
			& {Y\times_X \operatorname{Spec}(k^\prime)} &&&& {Y.}
			\arrow["{q^{\prime\prime}_2}", from=1-1, to=1-4]
			\arrow["{q^{\prime\prime}_1}", from=1-1, to=2-2]
			\arrow["{s_\ell}", from=1-1, to=3-1]
			\arrow["{s_1}", from=1-4, to=2-6]
			\arrow["{s_k}"{pos=0.8}, from=1-4, to=3-4]
			\arrow["{s_2}"{pos=0.2}, from=2-2, to=2-6]
			\arrow["{s_{k^\prime}}"{pos=0.8}, from=2-2, to=4-2]
			\arrow["s", from=2-6, to=4-6]
			\arrow["{q_2^\prime}"{pos=0.7}, from=3-1, to=3-4]
			\arrow["{q_1^\prime}"', from=3-1, to=4-2]
			\arrow["{g_k}", from=3-4, to=4-6]
			\arrow["{g_{k^\prime}}"', from=4-2, to=4-6]
		\end{tikzcd}
	\end{displaymath}
	Base change implies $s_k$, and hence $s_{\ell}$, are \'{e}tale presentations.
	Hence, $f_\ell \circ s_\ell$ is smooth and finitely presented.
	Denote by $T_\# $ the set of points of $|U\times_X \operatorname{Spec}(\# )|$ where $f_\#  \circ s_\# $ is smooth where $\#  \in \{ k^\prime, \ell \}$.
	Applying \cite[\href{https://stacks.math.columbia.edu/tag/02V4}{Tag 02V4}]{StacksProject}, $(q^{\prime\prime}_1)^{-1} (T_{k^\prime}) = T_\ell$.
	Since $q_1$ is surjective, so is $q^{\prime\prime}_1$.
	Hence, $q^{\prime\prime}_1 ((q^{\prime\prime}_1)^{-1} (T_{k^\prime})) = T_{k^\prime}$.
	However, $T_\ell = |U\times_X \operatorname{Spec}(\ell)|$, and so $f_{k^\prime} \circ s_{k^\prime}$ is smooth.
	As $s_{k^\prime}$ is an \'{e}tale presentation from a scheme, \cite[\href{https://stacks.math.columbia.edu/tag/03ZF}{Tag 03ZF}]{StacksProject} says $f_{k^\prime}$ is smooth, which completes the proof.
\end{proof}

\begin{lemma}
	\label{lem:smooth_characterize}
	Let $f\colon Y \to X$ be a finitely presented flat morphism of Noetherian algebraic spaces.
	Then the following are equivalent:
	\begin{enumerate}
		\item \label{lem:smooth_characterize1} $f$ is smooth
		\item \label{lem:smooth_characterize2} the base change of $f$ along every $\operatorname{Spec}(k)\to X$ is smooth where $k$ is a field
		\item \label{lem:smooth_characterize3} for every $p\in |X|$ there exists a representative $\operatorname{Spec}(k) \to X$ of $p$ such that the base change of $f$ along $\operatorname{Spec}(k)\to X$ is smooth.
	\end{enumerate}
\end{lemma}

\begin{proof}
	$\eqref{lem:smooth_characterize1} \implies \eqref{lem:smooth_characterize3}$ follows from \cite[\href{https://stacks.math.columbia.edu/tag/03ZF}{Tag 03ZF}]{StacksProject}, $\eqref{lem:smooth_characterize3} \implies \eqref{lem:smooth_characterize2}$ is \Cref{lem:smooth_if_some_representatives}, and $\eqref{lem:smooth_characterize2} \implies \eqref{lem:smooth_characterize1}$ is \Cref{lem:smooth_if_all_representatives}.
\end{proof}

\begin{lemma}
	\label{lem:regular_iff_for_all_field_extensions}
	Let $k$ be a field.
	A finitely presented morphism $f\colon X \to \operatorname{Spec}(k)$ of algebraic spaces is smooth if, and only if, $X\times_k \operatorname{Spec}(k^\prime)$ is regular for all field extensions $k^\prime/k$.
\end{lemma}

\begin{proof}
	If $f$ is smooth, then the base change of $f$ along all field extensions $k^\prime/k$ is smooth, and hence, $X\times_k \operatorname{Spec}(k^\prime)$ is regular.
	Conversely, let $s\colon U \to X$ be an \'{e}tale presentation.
	Fix a field extension $k^\prime/k$.
	If $X\times_k \operatorname{Spec}(k^\prime)$ is regular, then $U\times_k \operatorname{Spec}(k^\prime)$ is regular because the base change of $s$ along $k^\prime/k$ yields an \'{e}tale presentation.
	By \cite[\href{https://stacks.math.columbia.edu/tag/038V}{Tag 038V}]{StacksProject}, $U$ is geometrically regular.
	Then \cite[\href{https://stacks.math.columbia.edu/tag/038X}{Tag 038X}]{StacksProject} implies $U$ is smooth over $k$, i.e.\ $f \circ s$ is smooth.
	Thus, \cite[\href{https://stacks.math.columbia.edu/tag/03ZF}{Tag 03ZF}]{StacksProject} says $f$ is smooth.
\end{proof}

\begin{remark}
	\label{rmk:regular_locus_generalization}
	Let $X$ be a Noetherian algebraic space.
	Recall that the locus of $x\in |X|$ for which the \'{e}tale local ring $\mathcal{O}_{X,\overline{x}}$ is regular is stable under generalization.
	Indeed, choose an \'{e}tale presentation $s\colon U \to X$. 
	Let $x \rightsquigarrow x^\prime$ be a specialization in $|X|$ such that $\mathcal{O}_{X,\overline{x^\prime}}$ is regular.
	By \cite[\href{https://stacks.math.columbia.edu/tag/03JX}{Tag 03JX}, \href{https://stacks.math.columbia.edu/tag/03JV}{03JV}, and \href{https://stacks.math.columbia.edu/tag/03K2}{03K2}]{StacksProject}, there exists a specialization $u \rightsquigarrow u^\prime$ in $|U|$ such that $s(u) = x$ and $s(u^\prime) = x^\prime$.
	By \cite[\href{https://stacks.math.columbia.edu/tag/04KF}{Tag 04KF} and \href{https://stacks.math.columbia.edu/tag/06LN}{06LN}]{StacksProject}, $\mathcal{O}_{U,u^\prime}$ is regular.
	Hence, $\mathcal{O}_{U,u}$ is regular because it is a localization of $\mathcal{O}_{U,u^\prime}$ \cite[\href{https://stacks.math.columbia.edu/tag/0AFS}{Tag 0AFS}]{StacksProject}.
	Once more, from \cite[\href{https://stacks.math.columbia.edu/tag/04KF}{Tag 04KF} and \href{https://stacks.math.columbia.edu/tag/06LN}{06LN}]{StacksProject}, $\mathcal{O}_{X,\overline{x}}$ is regular.
\end{remark}

\begin{lemma}
	\label{lem:regular_iff_perfect_is_dbcoh}
	A Noetherian algebraic space $X$ is regular if, and only if, $\operatorname{Perf}(X) = D^b_{\operatorname{coh}}(X)$.
\end{lemma}

\begin{proof}
	Fix an \'{e}tale presentation $s\colon U \to X$ with $U$ affine.
	If $X$ is regular, then $U$ is regular.
	Hence, for any $E\in D^b_{\operatorname{coh}}(X)$, $\mathbf{L}s^\ast E$ is perfect \cite[\href{https://stacks.math.columbia.edu/tag/0FDC}{Tag 0FDC}]{StacksProject}, and so $E$ is perfect \cite[Lemma 3.5]{GuisadoVillaalgordo/Lank:2026}.
	It follows that $\operatorname{Perf}(X) = D^b_{\operatorname{coh}}(X)$.
	We prove the converse.
	Note that every point of $|X|$ specializes to a closed point \cite[\href{https://stacks.math.columbia.edu/tag/03K6}{Tags 03K6} and \href{https://stacks.math.columbia.edu/tag/005E}{005E}]{StacksProject}.
	Hence, by \Cref{rmk:regular_locus_generalization}, it suffices to prove that every closed point of $|X|$ belongs to the locus of points whose \'{e}tale local ring is regular.
	Let $p\in |X|$ be closed.
	Denote by $t_p\colon \operatorname{Spec}(\kappa(p))\to X$ the associated residual space of $p$ \cite[\href{https://stacks.math.columbia.edu/tag/06R0}{Tag 06R0}, \href{https://stacks.math.columbia.edu/tag/0H1R}{0H1R}, and \href{https://stacks.math.columbia.edu/tag/03I8}{03I8}]{StacksProject}.
	By \cite[\href{https://stacks.math.columbia.edu/tag/0H1U}{Tag 0H1U}]{StacksProject}, $t_p$ is a closed immersion.
	There exists a fibered square
	\begin{displaymath}
		\begin{tikzcd}
			{U\times_X \operatorname{Spec}(\kappa(p))} & {\operatorname{Spec}(\kappa(p))} \\
			U & {X.}
			\arrow["{s_p}", from=1-1, to=1-2]
			\arrow["{t^\prime_p}"', from=1-1, to=2-1]
			\arrow["{t_p}", from=1-2, to=2-2]
			\arrow["s"', from=2-1, to=2-2]
		\end{tikzcd}
	\end{displaymath}
	By \cite[Lemma 2.19]{GuisadoVillaalgordo/Lank:2026}, $s$ is quasi-affine, and hence representable.
	Base change implies $s_p$ is an \'{e}tale presentation. 
	By \cite[\href{https://stacks.math.columbia.edu/tag/02GL}{Tag 02GL}]{StacksProject}, $U\times_X \operatorname{Spec}(\kappa(p))$ is a disjoint union of spectra of finite separable field extensions of $\kappa(p)$.
	As $s_p$ is finitely presented, this disjoint union is finite, and hence its source is affine and regular.
	By \Cref{ex:regular_ring}, $\mathcal{O}_{U\times_X \operatorname{Spec}(\kappa(p))}$ finitely builds every object of $D^b_{\operatorname{coh}}(U\times_X \operatorname{Spec}(\kappa(p)))$.
	As $t_p$ is a closed immersion, $\mathbf{R}(t_p)_\ast \mathcal{O}_{\operatorname{Spec}(\kappa(p))}$ belongs to $D^b_{\operatorname{coh}}(X)$, which equals $\operatorname{Perf}(X)$ by hypothesis.
	By flat base change, $\mathbf{R}(t^\prime_p)_\ast \mathcal{O}_{U\times_X \operatorname{Spec}(\kappa(p))}$ is perfect.
	Hence, for all $E\in D^b_{\operatorname{coh}}(U\times_X \operatorname{Spec}(\kappa(p)))$, $\mathbf{R}(t^\prime_p)_\ast E$ is perfect. 
	Choose a closed point in $|U\times_X \operatorname{Spec}(\kappa(p))|$.
	Set $u\in |U|$ to be its image under $|t^\prime_p|$. 
	Denote by $h_u \colon \operatorname{Spec}(\mathcal{O}_{U,u}) \to U$ the canonical morphism.
	Consider the fibered square 
	\begin{displaymath}
		\begin{tikzcd}
			{\operatorname{Spec}(\mathcal{O}_{U,u})\times_X \operatorname{Spec}(\kappa(p))} & {U\times_X \operatorname{Spec}(\kappa(p))} \\
			{\operatorname{Spec}(\mathcal{O}_{U,u})} & {U.}
			\arrow["{h^\prime_u}", from=1-1, to=1-2]
			\arrow["{t^{\prime \prime}_p}"', from=1-1, to=2-1]
			\arrow["{t^\prime_p}", from=1-2, to=2-2]
			\arrow["{h_u}"', from=2-1, to=2-2]
		\end{tikzcd}
	\end{displaymath}
	Base change implies $t^{\prime \prime}_p$ is a closed immersion.
	By \cite[Lemma 4.8]{Lank:2026}, $\mathbf{L}(h^\prime_u)^\ast$ is essentially dense on $D^b_{\operatorname{coh}}$.
	Coupled with flat base change, it follows that 
	\begin{displaymath}
		\begin{aligned}
			\mathbf{R}(t^{\prime \prime}_p)_\ast D^b_{\operatorname{coh}}(\operatorname{Spec}(\mathcal{O}_{U,u})\times_X \operatorname{Spec}(\kappa(p))) 
			&\subseteq \langle \mathbf{R}(t^{\prime \prime}_p)_\ast  \mathbf{L}(h^\prime_u)^\ast D^b_{\operatorname{coh}}(U\times_X \operatorname{Spec}(\kappa(p))) \rangle
			\\&\subseteq \langle \mathbf{L}h_u^\ast \mathbf{R}(t^\prime_p)_\ast D^b_{\operatorname{coh}} (U\times_X \operatorname{Spec}(\kappa(p))) \rangle
			\\&\subseteq \langle \mathbf{L}h_u^\ast \operatorname{Perf}(U) \rangle 
			\\&\subseteq \operatorname{Perf}(\mathcal{O}_{U,u}).
		\end{aligned}
	\end{displaymath}
	Denote by $g_u \colon \operatorname{Spec}(\kappa(u)) \to \operatorname{Spec}(\mathcal{O}_{U,u})$ the canonical morphism.
	By \cite[\href{https://stacks.math.columbia.edu/tag/0356}{Tag 0356}]{StacksProject}, $g_u$ factors through $t^{\prime \prime}_p$.
	Since $g_u$ and $t^{\prime \prime}_p$ are closed immersions, the induced morphism from $\operatorname{Spec}(\kappa(u))$ to the source of $t^{\prime \prime}_p$ is a closed immersion.
	Its structure sheaf pushforward is coherent, and pushforwards along closed immersions are exact \cite[\href{https://stacks.math.columbia.edu/tag/04CJ}{Tag 04CJ}]{StacksProject}.
	Thus, $\mathbf{R}(g_u)_\ast \mathcal{O}_{\operatorname{Spec}(\kappa(u))}$ belongs to the left-hand side of the above inclusion, and hence is perfect.
	This implies $\mathcal{O}_{U,u}$ is a regular local ring \cite[\href{https://stacks.math.columbia.edu/tag/00OC}{Tag 00OC}]{StacksProject}.
	Hence, $u$ is a regular point of $U$.
	As $s(u) = p$, it follows that $X$ is regular at $p$ \cite[\href{https://stacks.math.columbia.edu/tag/0AH9}{Tag 0AH9}]{StacksProject}.
\end{proof}

\begin{proposition}
	\label{prop:smoothness}
	Let $f_i\colon Y_i\to S$ be proper flat morphisms of Noetherian algebraic spaces.
	Suppose $Y_1$ and $Y_2$ are Fourier--Mukai $S$-partners given by a kernel $E\in D^b_{\operatorname{coh}}(Y_1\times_S Y_2)$.
	Then $f_1$ is smooth if, and only if, $f_2$ is smooth.
\end{proposition}

\begin{proof}
	By \cite[Remark 2.7]{GuisadoVillaalgordo/Lank:2026}, $E$ is relatively perfect over each $Y_i$.
	By symmetry, it suffices to prove the case where $f_1$ is smooth.
	Fix a morphism $\operatorname{Spec}(k) \to S$ from a field.
	Let $L/k$ be a field extension.
	The composition $\operatorname{Spec}(L) \to S$ is affine \cite[\href{https://stacks.math.columbia.edu/tag/09TF}{Tag 09TF}]{StacksProject}.
	By \cite[Theorem 10.13]{GuisadoVillaalgordo/Lank:2026}, $Y_1 \times_S \operatorname{Spec}(L)$ and $Y_2 \times_S \operatorname{Spec}(L)$ are Fourier--Mukai $\operatorname{Spec}(L)$-partners.
	The kernel is the derived pullback of $E$ to $Y_1\times_S Y_2\times_S \operatorname{Spec}(L)$, which is relatively perfect over each $Y_i \times_S \operatorname{Spec}(L)$ \cite[Corollary 6.18]{GuisadoVillaalgordo/Lank:2026}.
	Consequently, \cite[Lemma 10.9]{GuisadoVillaalgordo/Lank:2026} implies the associated integral transform restricts to equivalences on both $D^b_{\operatorname{coh}}$ and $\operatorname{Perf}$.
	As $Y_1 \times_S \operatorname{Spec}(L)$ is smooth over $L$, it is regular.
	By \Cref{lem:regular_iff_perfect_is_dbcoh}, it follows that
	\begin{displaymath}
		\operatorname{Perf}(Y_1 \times_S \operatorname{Spec}(L))
		= D^b_{\operatorname{coh}} (Y_1 \times_S \operatorname{Spec}(L)).
	\end{displaymath}
	The integral transform takes the left-hand side to $\operatorname{Perf}(Y_2 \times_S \operatorname{Spec}(L))$ and the right-hand side to $D^b_{\operatorname{coh}}(Y_2 \times_S \operatorname{Spec}(L))$.
	Hence, every object of $D^b_{\operatorname{coh}} (Y_2 \times_S \operatorname{Spec}(L))$ is perfect, and so \Cref{lem:regular_iff_perfect_is_dbcoh} implies $Y_2 \times_S \operatorname{Spec}(L)$ is regular.
	Applying \Cref{lem:regular_iff_for_all_field_extensions}, $Y_2 \times_S \operatorname{Spec}(k)$ is smooth over $k$ because $L$ was arbitrary.
	Therefore, by \Cref{lem:smooth_characterize}, $f_2$ is smooth.
\end{proof}

\begin{lemma}
	\label{lem:obstruction_to_singular_equivalence}
	Let $f_1\colon Y_1 \to S$ and $f_2\colon Y_2 \to S$ be proper flat morphisms to a regular Noetherian algebraic space.
	Suppose $Y_1$ and $Y_2$ are Fourier--Mukai $S$-partners given by a kernel $B\in D^b_{\operatorname{coh}}(Y_1\times_S Y_2)$.
	If $Y_1$ is not regular, then $B\not\in\operatorname{Perf}(Y_1\times_S Y_2)$.
\end{lemma}

\begin{proof}
	By \cite[Remark 2.7]{GuisadoVillaalgordo/Lank:2026}, $B$ is relatively perfect over each $Y_i$.
	We prove the claim by contradiction.
	Assume that $B\in \operatorname{Perf}(Y_1\times_S Y_2)$.
	By \cite[Lemma 10.9]{GuisadoVillaalgordo/Lank:2026}, $\Phi_B$ restricts to triangulated equivalences on $D^b_{\operatorname{coh}}$ and $\operatorname{Perf}$.
	We first show that $\Phi_B (D^b_{\operatorname{coh}}(Y_1))\subseteq \operatorname{Perf}(Y_2)$.
	Choose a classical generator $G_i$ for $\operatorname{Perf}(Y_i)$, see \cite[\href{https://stacks.math.columbia.edu/tag/09IY}{Tags 09IY}, \href{https://stacks.math.columbia.edu/tag/09M8}{09M8}, and \href{https://stacks.math.columbia.edu/tag/09SR}{09SR}]{StacksProject}.
	By \cite[Lemma 5.4]{GuisadoVillaalgordo/Lank:2026} and \cite[\href{https://stacks.math.columbia.edu/tag/09SR}{Tag 09SR}]{StacksProject}, $\mathbf{L} p_1^\ast G_1 \otimes^{\mathbf{L}} \mathbf{L} p_2^\ast G_2$ is a classical generator for $\operatorname{Perf}(Y_1 \times_S Y_2)$.
	In particular, $B$ is finitely built by $\mathbf{L} p_1^\ast G_1 \otimes^{\mathbf{L}} \mathbf{L} p_2^\ast G_2$.
	Fix $E\in D^b_{\operatorname{coh}}(Y_1)$.
	As $\mathbf{R}(p_2)_\ast (\mathbf{L} p_1^\ast E \otimes^{\mathbf{L}} -)$ is an exact functor, it follows that
	\begin{displaymath}
		\Phi_B(E) \in \langle \mathbf{R}(p_2)_\ast (\mathbf{L} p_1^\ast E \otimes^{\mathbf{L}} \mathbf{L} p_1^\ast G_1 \otimes^{\mathbf{L}} \mathbf{L} p_2^\ast G_2) \rangle.
	\end{displaymath}
	By the projection formula and flat base change, there are isomorphisms
	\begin{displaymath}
		\begin{aligned}
			\mathbf{R}(p_2)_\ast (\mathbf{L} p_1^\ast E \otimes^{\mathbf{L}} \mathbf{L} p_1^\ast G_1 \otimes^{\mathbf{L}} \mathbf{L} p_2^\ast G_2)
			&\cong \mathbf{R}(p_2)_\ast \mathbf{L} p_1^\ast (E \otimes^{\mathbf{L}} G_1) \otimes^{\mathbf{L}} G_2
			\\&\cong \mathbf{L} f_2^\ast \mathbf{R}(f_1)_\ast (E \otimes^{\mathbf{L}} G_1) \otimes^{\mathbf{L}} G_2.
		\end{aligned}
	\end{displaymath}
	As $f_1$ is proper, $\mathbf{R}(f_1)_\ast (E \otimes^{\mathbf{L}} G_1)\in D^b_{\operatorname{coh}}(S)$.
	Since $S$ is regular, \Cref{lem:regular_iff_perfect_is_dbcoh} implies $D^b_{\operatorname{coh}}(S)= \operatorname{Perf}(S)$.
	Hence, $\Phi_B(E)$ is perfect, as desired.
	Now, $\Phi_B$ takes $D^b_{\operatorname{coh}}(Y_1)$ to $D^b_{\operatorname{coh}}(Y_2)$, so $D^b_{\operatorname{coh}}(Y_2) = \operatorname{Perf}(Y_2)$.
	As $\Phi_B$ restricts to equivalences on $D^b_{\operatorname{coh}}$ and $\operatorname{Perf}$, it follows that $D^b_{\operatorname{coh}}(Y_1) = \operatorname{Perf}(Y_1)$.
	Thus, \Cref{lem:regular_iff_perfect_is_dbcoh} implies $Y_1$ is regular, which is a contradiction.
\end{proof}

\begin{lemma}
	\label{lem:perf_by_fibers_of_morphism}
	Let $f\colon Y \to X$ be a finitely presented morphism of Noetherian algebraic spaces.
	For each $p\in |X|$, denote by $t_p\colon \operatorname{Spec}(\kappa(p)) \to X$ the canonical morphism from the residue field of $p$, and by $t^\prime_p\colon Y\times_X \operatorname{Spec}(\kappa(p)) \to Y$ the morphism obtained by base change along $t_p$.
	If $E\in D^b_{\operatorname{coh}}(Y)$ satisfies $\mathbf{L}(t^\prime_p)^\ast E\in \operatorname{Perf}(Y\times_X \operatorname{Spec}(\kappa(p)))$ for all $p\in |X|$, then $E\in \operatorname{Perf}(Y)$.
\end{lemma}

\begin{proof}
	As $f$ is finitely presented, $Y\times_X \operatorname{Spec}(\kappa(p))$ is Noetherian for all $p\in |X|$.
	Choose $q\in |Y|$, and set $p := f(q)$.
	By \cite[\href{https://stacks.math.columbia.edu/tag/03H4}{Tag 03H4}]{StacksProject}, there exists $q^\prime\in |Y\times_X \operatorname{Spec}(\kappa(p))|$ which maps to $q$.
	Denote by $b\colon \operatorname{Spec}(\kappa(q^\prime)) \to Y\times_X \operatorname{Spec}(\kappa(p))$ the canonical morphism from the residue field of $q^\prime$, see \cite[\href{https://stacks.math.columbia.edu/tag/0EMV}{Tag 0EMV}]{StacksProject}.
	The image of $t^\prime_p\circ b$ is $q$, and so $t^\prime_p \circ b$ represents $q$ \cite[Lemma 2.8]{GuisadoVillaalgordo/Lank:2026}.
	As $\mathbf{L}(t^\prime_p)^\ast E$ is perfect, so is $\mathbf{L}(t^\prime_p \circ b)^\ast E \cong \mathbf{L}b^\ast \mathbf{L}(t^\prime_p)^\ast E$.
	In particular, $\mathbf{L}(t^\prime_p \circ b)^\ast E$ has bounded cohomology.
	Since $q$ was arbitrary, \cite[Proposition 3.12]{GuisadoVillaalgordo/Lank:2026} implies $E$ is perfect.
\end{proof}

\begin{proposition}
	\label{prop:kernel_smoothness}
	Let $f_i\colon Y_i\to S$ be proper flat morphisms of Noetherian algebraic spaces.
	Suppose $Y_1$ and $Y_2$ are Fourier--Mukai $S$-partners given by a kernel $E\in D^b_{\operatorname{coh}}(Y_1\times_S Y_2)$.
	Then $f_1$ is smooth if, and only if, $E\in\operatorname{Perf}(Y_1\times_S Y_2)$.
\end{proposition}

\begin{proof}
	By \cite[Remark 2.7]{GuisadoVillaalgordo/Lank:2026}, $E$ is relatively perfect over each $Y_i$.

	Assume that $f_1$ is smooth.
	Choose $p\in |S|$.
	Denote by $t_p\colon \operatorname{Spec}(\kappa(p)) \to S$ the canonical morphism associated to the residue field of $p$.
	Denote by $t^\prime_p\colon Y_1\times_S Y_2 \times_S \operatorname{Spec}(\kappa(p)) \to Y_1\times_S Y_2$ the canonical morphism obtained by base change along $t_p$.
	As $f_1$ is smooth, \Cref{prop:smoothness} implies $f_2$ is smooth.
	Hence, by base change, $Y_1\times_S Y_2 \times_S \operatorname{Spec}(\kappa(p))$ is smooth over $\kappa(p)$.
	Then \Cref{lem:regular_iff_perfect_is_dbcoh} implies
	\begin{displaymath}
		\operatorname{Perf}(Y_1\times_S Y_2 \times_S \operatorname{Spec}(\kappa(p))) = D^b_{\operatorname{coh}}(Y_1\times_S Y_2 \times_S \operatorname{Spec}(\kappa(p))).
	\end{displaymath}
	By \cite[Corollary 6.18]{GuisadoVillaalgordo/Lank:2026}, $\mathbf{L}(t^\prime_p)^\ast E$ is relatively perfect over each $Y_i \times_S \operatorname{Spec}(\kappa(p))$.
	In particular, $\mathbf{L}(t^\prime_p)^\ast E$ has bounded cohomology, and so it is perfect.
	Since this holds for all $p\in |S|$, \Cref{lem:perf_by_fibers_of_morphism} asserts that $E$ is perfect.

	Assume that $E$ is perfect.
	Choose $p\in |S|$.
	Denote by $t_p\colon \operatorname{Spec}(\overline{\kappa(p)}) \to S$ a  morphism from an algebraic closure of the residue field of $p$.
	Denote by $t^\prime_p\colon Y_1\times_S Y_2 \times_S \operatorname{Spec}(\overline{\kappa(p)}) \to Y_1\times_S Y_2$ the canonical morphism obtained by base change along $t_p$.
	By \cite[Corollary 6.18]{GuisadoVillaalgordo/Lank:2026}, $\mathbf{L}(t^\prime_p)^\ast E$ is relatively perfect over each $Y_i \times_S \operatorname{Spec}(\overline{\kappa(p)})$.
	Moreover, $t_p$ is affine \cite[\href{https://stacks.math.columbia.edu/tag/09TF}{Tag 09TF}]{StacksProject}, so \cite[Theorem 10.13]{GuisadoVillaalgordo/Lank:2026} implies $\mathbf{L}(t^\prime_p)^\ast E$ induces an equivalence on $D^b_{\operatorname{coh}}$.
	Furthermore, $\mathbf{L}(t^\prime_p)^\ast E$ is perfect, and so \Cref{lem:obstruction_to_singular_equivalence} implies $Y_1 \times_S \operatorname{Spec}(\overline{\kappa(p)})$ is regular.
	As $\overline{\kappa(p)}$ is perfect, it has no nontrivial finite purely inseparable extensions.
	Applying \cite[\href{https://stacks.math.columbia.edu/tag/038V}{Tags 038V} and \href{https://stacks.math.columbia.edu/tag/038X}{038X}]{StacksProject} to an \'{e}tale presentation of $Y_1 \times_S \operatorname{Spec}(\overline{\kappa(p)})$, it follows from \cite[\href{https://stacks.math.columbia.edu/tag/03ZF}{Tag 03ZF}]{StacksProject} that $Y_1 \times_S \operatorname{Spec}(\overline{\kappa(p)})$ is smooth over $\overline{\kappa(p)}$.
	Since $t_p$ represents $p$ and $p\in |S|$ was arbitrary, \Cref{lem:smooth_characterize} asserts that $f_1$ is smooth, which completes the proof.
\end{proof}

\bibliographystyle{alpha}
\bibliography{mainbib}

\end{document}